\documentclass[]{article}
\usepackage{epic,eepic,graphicx,color}        
\usepackage[bottom]{footmisc}
\usepackage{amsmath,amssymb,amsthm,bm}

\usepackage{multirow}

\usepackage[left=1cm,right=1cm,top=1cm,bottom=1cm]{geometry}

\begin{document}

\title{Exact pressure elimination for the Crouzeix-Raviart scheme\\ applied to the Stokes and Navier-Stokes problems}
\author{Eric Ch\'enier\thanks{MSME, Univ. Gustave Eiffel, Univ. Paris Est Créteil, CNRS,
F-77454, Marne-la-Vall\'ee, France. \tt{eric.chenier@univ-eiffel.fr}}~ and
Robert Eymard\thanks{ LAMA, Univ. Gustave Eiffel, Univ. Paris Est Créteil, CNRS,
F-77454, Marne-la-Vall\'ee, France. \tt{Robert.Eymard@univ-eiffel.fr}}}

\maketitle

\begin{abstract} We show that, using the Crouzeix-Raviart scheme, a cheap algebraic transformation, applied to the coupled velocity--pressure linear systems issued from the transient or steady Stokes or Navier-Stokes problems, leads to a linear system only involving as many auxiliary variables as the velocity components. This linear system, which is symmetric positive definite in the case of the transient Stokes problem and symmetric invertible in the case of the steady Stokes problem, with the same stencil as that of the velocity matrix, provides the exact solution of the initial coupled linear system.
Numerical results show the increase of performance when applying direct or iterative solvers to the resolution of these linear systems.
\end{abstract}

\section{Introduction}

This paper is focused on the resolution of the coupled velocity-pressure linear systems issued from the discretisation by the Crouzeix-Raviart scheme \cite{crou-73-con} of the steady (or  transient and semi-discretised in time) Stokes and Navier-Stokes problems, considering for the sake of simplicity homogeneous Dirichlet boundary conditions for the velocity. In the case of the 
 Navier-Stokes problem, these linear systems are resulting from the Newton-Raphson method applied to iteratively solve the non-linear equations. 

\medskip

These linear systems are under the form
\begin{equation}\label{eq:syslinintro}
\begin{bmatrix} A & D^t \\ D & 0\end{bmatrix} \begin{bmatrix} U\\ P\end{bmatrix}=  \begin{bmatrix} R\\ 0\end{bmatrix},
\end{equation}
where $A$ is the rigidity matrix resulting from the use of the ${\mathbb P}^1$ non conforming finite element method for the velocities, completed by the mass matrix in case of the transient case, and by some derivatives issued from the convection term in case of the Navier-Stokes problem, $D$ is the discrete divergence matrix written element by element, $U$ is the vector of all velocity unknowns, $P$ is the vector of all (but one) pressures unknowns and $R$ is the right hand-side resulting from the momentum source terms.

\medskip

In the case of the steady or transient Stokes problem, the matrix $A$ is symmetric. 
This is no longer the case for the Navier-Stokes problem. 
But even in the case of the  steady or transient Stokes problem, the matrix of the linear system \eqref{eq:syslinintro} is not positive definite, due to the fact that there are negative eigenvalues since there are zeros on the main diagonal. 
This property makes much more complicate the use of iterative solvers based, for example, on conjugate gradient or  GMRES \cite{gmres} methods preconditioned with Incomplete Lower-Upper (ILU) factorization. 
Note that the implementation of the ILU preconditioners on parallel architectures \cite{cho2015fine} fails to provide the same preconditioning properties as ILU on only one processor, due to the loss of some sequential computations.

\medskip

Then many authors are led, on small cases, to use direct solvers (recall the remarkable performances of the direct MUMPS solvers on parallel architectures \cite{mumps,mumps2}). But on large matrices, such direct methods can no longer be reasonably applied, and there is a need to use all the same an iterative linear solver.

\medskip

Another option consists in adding a small diagonal pressure-pressure connection, as performed by the augmented Lagrangian methods.
But then, the iterative convergence properties of the solutions for such a modified system to that of the original one may become very slow.


\medskip

Such difficulties for solving the linear systems issued from a mixed formulation are well-known when solving a simple Laplace problem. In this case, $H_{\rm div}$ conforming finite elements are used for the approximation of the gradient of the unknown (the Raviart-Thomas finite element is often used  in the case of simplicial meshes), and piecewise constant elements are used for the unknown. A very clever method is then known for overcoming the difficulty of solving the linear systems issued from this problem: it is the famous hybridisation of the problem, leading to solve a symmetric positive definite linear system on the trace of the unknown on the faces of the mesh \cite{chen1996equi, duran, vohralik}. Note that a similar idea is used in \cite{agh2015hyb} in the case of the Stokes problem, discretised by Hybrid High Order methods.

\medskip

This paper is based on the extension of the same idea for applying an algebraic hybridisation to the case of the coupled linear systems \eqref{eq:syslinintro}. Let us emphasize that the solution of the linear system is not modified by the use of this hybridisation. In order the method to apply to the Navier-Stokes problem, we select an implementation of the non-linear convection term which does not increase the stencil of the Stokes problem \cite{herlat}. 

\medskip

In the transient Stokes problem, we get, after hybridisation, a symmetric positive definite linear system with as many unknowns as the velocities, and the same connection stencil (even in the case of the Stokes problem, the different space components of the auxiliary unknowns are connected, contrarily to the original velocity-velocity matrix).

\medskip

In the steady case, we are led to introduce a modification in the diagonal blocks  to have an invertible block diagonal matrix. Once again, the solution of the linear system is not altered by this modification. After hybridisation, we obtain in the case of the Stokes problem final symmetric linear system to be solved with as many unknowns as the velocities, but the matrix is no longer positive definite.

\medskip

This paper is organised as follows. We first detail in Section \ref{sec:crscheme} the construction of the scheme, with precising the treatment of the right-hand-side allowing exact numerical solutions in the case where it resumes to the gradient of a scalar field, and with a formulation of the convection term which does not enlarge the stencil. We then show in section \ref{sec:matrix} how the linear systems issued from this scheme can be algebraically handled for obtaining smaller linear systems with the same sparsity. We finally compare, in  Section \ref{sec:num},  the numerical efficiency of different linear solvers, applied to the initial coupled linear system and applied to their algebraic transformation.

\section{The Crouzeix-Raviart scheme for $d=2$ or $d=3$}\label{sec:crscheme}

Let us first give the strong formulation of the Stokes and Navier-Stokes equations in their steady or semi-discrete transient versions:

\begin{equation}\left\{\begin{array}{rll}\label{pforss}
\mu  \overline{\bm{u}} -\nu \Delta\overline{\bm{u}}+\nabla \overline{p} + {\bm{b}}(\overline{\bm{u}}) &=\overline{\bm{f}}&\hbox{ in }  \Omega\\
{\rm div} \overline{\bm{u}} &=0 &\hbox{ in } \Omega\\
\overline{\bm{u}} &=0 &\hbox{ on } \partial\Omega\\
\int_\Omega \overline{p}({\bm x}){\rm d}{\bm x} = 0
\end{array}\right.\end{equation}
where $\overline{\bm{u}} = (\overline{u}^{(i)})_{i=1,\ldots,d}$ with $d=2$ or $d=3$ represents the velocity field, $\Delta\overline{\bm{u}} = (\Delta\overline{u}^{(i)})_{i=1,\ldots,d}$, $\overline{p}$ is the pressure, the
domain $\Omega$ with boundary $\partial\Omega$ is a bounded open set in ${\mathbb R}^d$, $\nu>0$ is the invert of the Reynolds number,  $\overline{\bm{f}}=(\overline{f}^{(i)})_{i=1,\ldots,d}$ is a given function defined on $\Omega$, $\nabla \overline{p} = \big(\partial_i \overline{p}\big)_{i=1,\ldots,d}$, ${\rm div} \overline{\bm{u}} = \sum_{i=1}^d\partial_i \overline{u}^{(i)}$.

\medskip

For the steady problem, $\mu = 0$ and in the case where the problem is transient, $\mu>0$ is the invert of the time step: then $\overline{\bm{f}}$ includes a term issued from the velocity at the beginning of the time step (and the transient problem is semi-discretised in time).

\medskip

For the transient or steady Stokes problems, we let
\begin{equation}\label{eq:noconv}
 {\bm{b}}(\overline{\bm{u}}) = 0,
\end{equation}
and for the Navier-Stokes problem, we define the non-linear convection term by
\begin{equation}\label{eq:conv}
 {\bm{b}}(\overline{\bm{u}}) = (\overline{\bm{u}}\cdot\nabla)\overline{\bm{u}} = \big(\sum_{j=1}^d \overline{u}^{(j)}\partial_j \overline{u}^{(i)}\big)_{i=1,\ldots,d}~.
\end{equation}
The standard weak formulation of Problem \eqref{pforss} is the following mixed one. Defining  $L^2_0(\Omega)$  as the set of elements of $L^2(\Omega)$ with null mean value on $\Omega$, this formulation is given by

\begin{equation}\left\{\begin{array}{rll}\label{eq:pweak}
\hbox{Find }&\overline{\bm{u}}\in H^1_0(\Omega)^d\hbox{ and }\overline{p}\in L^2_0(\Omega)\hbox{ such that }\\
\forall \overline{\bm{v}}\in  H^1_0(\Omega)^d,&\displaystyle \int_\Omega \Big(\mu  \overline{\bm{u}}\cdot \overline{\bm{v}} +\nu \nabla\overline{\bm{u}}:\nabla\overline{\bm{v}}-\overline{p}{\rm div} \overline{\bm{v}} + {\bm{b}}(\overline{\bm{u}})\cdot \overline{\bm{v}}\Big){\rm d}{\bm x} &=\displaystyle \int_\Omega\overline{\bm{f}}\cdot\overline{\bm{v}}{\rm d}{\bm x}\\
\forall\overline{q}\in L^2_0(\Omega),&\displaystyle \int_\Omega{\rm div} \overline{\bm{u}}\ \overline{q}{\rm d}{\bm x}&=0
\end{array}\right.\end{equation}

The Crouzeix-Raviart scheme \cite{crou-73-con} is the translation of the weak formulation \eqref{eq:pweak}  into discrete sets and operators applying on simplicial meshes (triangles in 2D, tetrahedra in 3D). It reads 
\begin{equation}\left\{\begin{array}{rll}\label{eq:pweakdis}
\hbox{Find }&{\bm{u}}\in (V_h)^d\hbox{ and }{p}\in Q_{h,0}\hbox{ such that }\\
\forall {\bm{v}}\in  (V_h)^d,&\displaystyle \int_\Omega \Big(\mu   \Pi_h {\bm{u}}\cdot  \Pi_h {\bm{v}} +\nu \nabla_h{\bm{u}}:\nabla_h{\bm{v}}-{p}{\rm div}_h {\bm{v}}\Big){\rm d}{\bm x}  + b_h({\bm{u}}, {\bm{v}}  )&=\displaystyle \int_\Omega\overline{\bm{f}}\cdot \widehat{\Pi}_h {\bm{v}}{\rm d}{\bm x}\\
\forall{q}\in Q_{h,0},&\displaystyle \int_\Omega{\rm div}_h {\bm{u}}\ {q}{\rm d}{\bm x}&=0.
\end{array}\right.\end{equation}
Let us define each of the discrete objects involved in \eqref{eq:pweakdis}.

\begin{enumerate}
 \item {\bf The finite dimensional space $V_h$.} \\
Let ${\mathcal M}$ be a simplicial mesh, that is a finite set of disjoint open simplices whose closure recovers $\Omega$. For $K\in {\mathcal M}$, we denote by ${\bm x}_K$ the centre of gravity of $K$. Denote by ${\mathcal F}$ the set of all faces (edges in 2D) of the mesh, that is partitioned into ${\mathcal F}_{\rm int}\cup{\mathcal F}_{\rm ext} $ (the set of interior and exterior faces), and denote for any $K\in {\mathcal M}$ by ${\mathcal F}_K$ the set of the faces of $K$. We denote by ${\mathcal F}_{K,\rm int} = {\mathcal F}_K\cap{\mathcal F}_{\rm int}$.

For any $\sigma\in{\mathcal F}_K$, we denote by ${\bm n}_{K,\sigma}$ the unit vector, normal to $\sigma$ and outward to $K$, and we let
\[
 {\bm a}_{K,\sigma} = |\sigma|{\bm n}_{K,\sigma}.
\]

We assume that there are no hanging nodes, which implies that the cardinal of any ${\mathcal F}_K$ is equal to $d+1$ (3 in 2D, 4 in 3D). For any face $\sigma\in {\mathcal F}$, we denote by ${\mathcal M}_\sigma$ the set of the simplices $K\in {\mathcal M}$ such that $\sigma\in {\mathcal F}_K$. Then the cardinal of ${\mathcal M}_\sigma$ is 2 for an interior face, 1 for an exterior face.
For any $\sigma\in{\mathcal F}$,  we denote by ${\bm x}_\sigma$ the centre of gravity of $\sigma$.

\medskip 

We then define, for any $\sigma\in {\mathcal F}_{\rm int}$ with ${\mathcal M}_\sigma = \{K,L\}$, the function $\varphi_{\sigma}~:~\Omega\to\mathbb{R}$ whose the restriction $\varphi_{\sigma,K}$ on $K$ (respectively $\varphi_{\sigma,L}$ on $L$)  is an affine function on $K$ (respectively $L$) and which is null on any other element of the mesh. Moreover, one requests that the mean values of both $\varphi_{\sigma,K}$ and $\varphi_{\sigma,L}$ are equal to $1$ on $\sigma$ and equal to $0$ on any $\sigma'\in {\mathcal F}_K\cup {\mathcal F}_L$ different from $\sigma$. These conditions are sufficient for defining in an unique way the affine functions $\varphi_{\sigma,K}$ and $\varphi_{\sigma,L}$, on each of which $d+1$ independent conditions have been specified. This definition ensures the continuity of the mean value of these functions on any face of the mesh, as well as the continuity of these functions at the centre of gravity of the faces of the mesh.

Then the space $V_h$ is defined as the space spanned by the family $(\varphi_{\sigma})_{\sigma\in {\mathcal F}_{\rm int}}$.

For any $v\in V_h$ and $K\in {\mathcal M}$, we denote by $v_K$ the restriction of $v$ to $K$ (it is therefore an affine function).

\medskip

 For any ${\bm v}\in (V_h)^d$ and  $\sigma\in {\mathcal F}$, we then denote by ${\bm V}_\sigma$ the vector ${\bm V}_\sigma  = (V_{\sigma}^{(i)} := v_K^{(i)}({\bm x}_\sigma))_{i=1,\ldots,d}$.

 \item {\bf The discrete operators $\nabla_h$ and ${\rm div}_h$.}\\
 The discrete operators $\nabla_h$ and ${\rm div}_h$ are defined as the ``broken'' ones, that means that there restriction to any element of the mesh are defined as the continuous ones:
\[
 \forall v\in {V}_{h},\ \forall K \in{\mathcal M}, \  (\nabla_h v)_{|K} = \nabla v_K\hbox{ and }\forall {\bm v}\in ({V}_{h})^d,\ \forall K \in{\mathcal M}, \  ({\rm div}_h {\bm v})_{|K} ={\rm div} v_K.
\]

\item  {\bf The discrete reconstruction operator $\Pi_h$.}\\
The operator $\Pi_h$ is introduced in order to obtain some mass lumping in the ``mass matrix'' term, that is in order to get a diagonal mass matrix. If $d=2$, using the Crouzeix-Raviart basis functions, the matrix 
\[
 M_{\sigma,\sigma'} = \int_K \varphi_\sigma({\bm x}) \varphi_{\sigma'}({\bm x}){\rm d}{\bm x}
\]
is already diagonal, and then $\Pi_h  = {\rm Id}$. But this fails if $d=3$. We then denote by $\Pi_h \varphi_\sigma$ a piecewise constant function, equal to 1 in a domain surrounding $\sigma$ and 0 elsewhere (this domain is defined as the union of the two triangles (2D) or tetrahedra (3D), the basis of which is $\sigma$, and the vertex of which is the centre of gravity of the neighbouring simplices).

\item  {\bf The discrete reconstruction operator $\widehat{\Pi}_h$.}\\
Following \cite{LIN14}, the operator $ \widehat{\Pi}_h {\bm{v}}$ is designed to ensure the following properties: $\widehat{\Pi}_h {\bm{v}}\in H_{\rm div}(\Omega)$ (which means a kind of continuity of the normal trace on any internal boundary), $\widehat{\Pi}_h {\bm{v}} -  {\bm{v}}$ tends to 0 as $h$ tends to 0 if ${\bm{v}}$ is the interpolation of any regular function, and finally there holds,
\begin{equation}\label{eq:linke}
 \forall \bm{v} \in  V_h^d,\ \left(\forall{q}\in Q_{h,0},\ \int_\Omega{\rm div}_h {\bm{v}} ~ {q}{\rm d}{\bm x}= 0\right)\Rightarrow {\rm div}\widehat{\Pi}_h {\bm{v}} = 0 \hbox{ a.e. in }\Omega.
\end{equation}

Indeed, if we change $\overline{\bm{f}}$ into $\overline{\bm{f}}+\nabla\varphi$, Property \eqref{eq:linke} implies that the discrete velocity is not modified, only the pressure field is changed by the addition of an interpolation of $\varphi$. 
This property leads to a substantial decrease of the numerical error, in particular in the case where the major part of $\overline{\bm{f}}$ is constituted by the gradient of a scalar field.
To this purpose, we use the Raviart-Thomas basis, which is conforming in $H_{\rm div}(\Omega)$ and  defined, for all $K\in {\mathcal M}$, $\sigma\in {\mathcal F}_K$ and ${\bm x}\in K$, by
\[
 {\bm \psi}_{K,\sigma}({\bm x}) = \frac {|\sigma|} {d\ |K|} ({\bm x} - {\bm s}_\sigma),
\]
where ${\bm s}_\sigma$ is the vertex of $K$ which is not a vertex of $\sigma$. Then we define, for any ${\bm x}\in K$,
\[
 \widehat{\Pi}_h {\bm v}({\bm x}) = \sum_{\sigma\in  {\mathcal F}_K} {\bm V}_\sigma\cdot {\bm n}_{K,\sigma} {\bm \psi}_{K,\sigma}({\bm x}) =  \sum_{\sigma\in  {\mathcal F}_K}  \frac {{\bm V}_\sigma\cdot {\bm a}_{K,\sigma}} {d\ |K|} ({\bm x} - {\bm s}_\sigma).
\]
We then approximate  $\int_K\overline{\bm{f}}\cdot \widehat{\Pi}_h {\bm{v}}{\rm d}{\bm x}$ by 
\begin{equation}
\int_\Omega\overline{\bm{f}}\cdot \widehat{\Pi}_h {\bm{v}}{\rm d}{\bm x} =  \sum_{K\in {\mathcal M}} \sum_{\sigma\in  {\mathcal F}_K}   {\bm V}_\sigma\cdot{\bm a}_{K,\sigma} \Big(\frac 1 {|K|}\int_K\overline{\bm f}{\rm d}{\bm x}\Big)\cdot ({\bm x}_\sigma - {\bm x}_K).
\label{eq:defrhs}\end{equation}
\item {\bf The finite dimensional space $Q_{h,0}$.}\\
We define $Q_h$ as the finite dimensional subset of $L^2(\Omega)$ spanned by the characteristic functions $\psi_K$ of all the simplices $K\in {\mathcal M}$ ($\psi_K$ is  the piecewise constant function defined on $\Omega$ which is equal to one inside $K$ and 0 elsewhere). 
Since the pressures can be defined up to a constant value, instead of defining a space of functions with null average (which would connect all components of the function together), we select a given element of the mesh ${\mathcal M}$, denoted $K_{0}$, and we define the set $Q_{h,0}$ as the set of all elements $p\in Q_h$ vanishing on $K_0$. Note that, for any $p\in Q_{h,0}$, we retrieve an element of  $L^2_0(\Omega)$, considering $p - \frac 1 {|\Omega|}\int_\Omega p({\bm x}){\rm d}{\bm x}$.

\item {\bf The non-linear form $b_h({\bm{u}}, {\bm{v}})$.}\\
This non-linear form vanishes for the transient or steady Stokes problems. For the Navier-Stokes problem, 
the following discretisation for $b_h({\bm{u}}, {\bm{v}}  )$ has been proposed by \cite{herlat} and is compared to other choices in \cite{efg}. Its main advantage is to keep a reduced stencil in the linear systems.
 All the simplices $K$ are split into co-volumes linked to the faces, as shown by Figure \ref{fig:covolume}. 
\begin{figure}[ht!]
\begin{center}
\resizebox{.4\textwidth}{!}{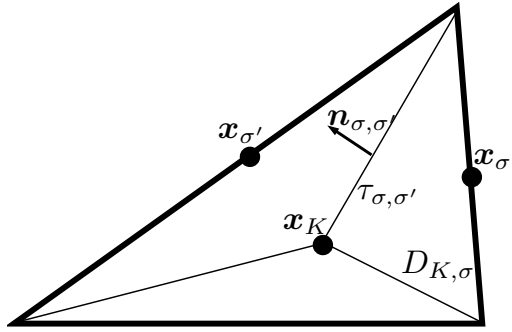}
\end{center}
\caption{Co-volumes associated with faces}\protect{\label{fig:covolume}}
\end{figure}
The co-volume associated with a face $\sigma$ in a simplex $K$, is defined as the cone $D_{K,\sigma}$ based on $\sigma$, whose vertex is the centre of gravity of $K$ (it is then a simplex as well). This sub-mesh leads to the definition of $d(d-1)$ internal faces, each of them being common to $D_{K,\sigma}$ and $D_{K,\sigma'}$, denoted $\tau_{\sigma,\sigma'}$, for any pair $\sigma,{\sigma'}\in {\mathcal F}_K$. Then the unit normal vector to the face  $\tau_{\sigma,\sigma'}$, oriented from $D_{K,\sigma}$ to $D_{K,\sigma'}$, is denoted by ${\bm n}_{\sigma,{\sigma'}}$. We then define  $b_h({\bm u},{\bm v})$ by the relation
\begin{equation}\label{eq:defbh}
b_h({\bm u},{\bm v}) := \sum_{K\in {\mathcal M}} \sum_{ \{ \sigma,\sigma' \}\subset {\mathcal F}_K }
 F_{\sigma,{\sigma'}}({\bm u}) ({\bm U}_{\sigma'}- {\bm U}_{{\sigma}})  \cdot \frac {{\bm V}_{\sigma}+ {\bm V}_{{\sigma'}}} {2},
\end{equation}
which also satisfies
\[
 b_h({\bm u},{\bm v}) = \frac 1 2 \sum_{K\in {\mathcal M}} \sum_{\sigma\in {\mathcal F}_K} {\bm V}_{\sigma} \cdot\sum_{{\sigma'}\in {\mathcal F}_K\setminus\{\sigma\}}
 F_{\sigma,{\sigma'}}({\bm u}) ({\bm U}_{\sigma'}- {\bm U}_{{\sigma}}),
\]

where $F_{\sigma,{\sigma'}}({\bm u})$ is defined by
\[
   F_{\sigma,{\sigma'}} ({\bm u}) = \int_{\tau_{\sigma,{\sigma'}}} {\bm u}_K ({\bm x})\cdot {\bm n}_{\sigma,{\sigma'}} {\rm d} s({\bm x}).
\]
We remark that, for $\sigma,{\sigma'}\in {\mathcal F}_K$, the centre of gravity ${\bm x}_{\sigma,{\sigma'}}$ of $\tau_{\sigma,{\sigma'}}$ is given by
\[
 {\bm x}_{\sigma,{\sigma'}} = {\bm x}_\sigma + {\bm x}_{\sigma'} - {\bm x}_K,
\]
and we observe that
\[
 |\tau_{\sigma,{\sigma'}}|{\bm n}_{\sigma,{\sigma'}} = \frac {1} {d+1}( {\bm a}_{K,\sigma'} - {\bm a}_{K,\sigma}).
\]
This yields
\[
  F_{\sigma,{\sigma'}} ({\bm u}) = \Big({\bm U}_{\sigma} + {\bm U}_{\sigma'} -  \frac {1} {d+1}\sum_{\sigma''\in  {\mathcal F}_K} {\bm U}_{\sigma''}\Big)\cdot\frac {1} {d+1}( {\bm a}_{K,\sigma'} - {\bm a}_{K,\sigma}).
\]
We then check that the relation $\sum_{\sigma'\in  {\mathcal F}_K} {\bm U}_{\sigma'}\cdot{\bm a}_{K,\sigma'} = 0$ implies that
\[
 \sum_{\sigma'\in  {\mathcal F}_K} F_{\sigma,{\sigma'}} ({\bm u}) =  -{\bm U}_{\sigma}\cdot{\bm a}_{K,\sigma}.
\]
Hence, the above definition is such that, if ${\rm div}_h {\bm u} = 0$, then there holds $b_h({\bm u},{\bm u})=0$ for all ${\bm u}\in V_h$. Indeed, there holds
\[
0 = \int_{D_{K,\sigma}} {\rm div} {\bm u}_K ({\bm x}) {\rm d} {\bm x} =\sum_{\sigma'\in {\mathcal F}_K\setminus\{\sigma\}}  F_{\sigma,{\sigma'}} ({\bm u}) + \int_{\sigma} {\bm u}_K ({\bm x})\cdot {\bm n}_{K,\sigma} {\rm d} s({\bm x}),
\]
which implies that
\begin{multline*}
 b_h({\bm u},{\bm u})= \frac 1 2\sum_{K\in {\mathcal M}}\sum_{ \{ \sigma,\sigma' \}\subset {\mathcal F}_K }
 F_{\sigma,{\sigma'}}({\bm u}) (|{\bm U}_{\sigma'}|^2- |{\bm U}_{{\sigma}}|^2)\\
 = -\frac 1 2 \sum_{K\in {\mathcal M}} \sum_{\sigma\in {\mathcal F}_K}|{\bm U}_{\sigma}|^2\sum_{\sigma'\in {\mathcal F}_K\setminus\{\sigma\}}  F_{\sigma,{\sigma'}} ({\bm u})
 =  \frac 1 2\sum_{K\in {\mathcal M}} \sum_{\sigma\in {\mathcal F}_K}|{\bm U}_{\sigma}|^2\int_{\sigma} {\bm u}_K ({\bm x})\cdot {\bm n}_{K,\sigma} {\rm d} s({\bm x}),
\end{multline*}
and this last term vanishes, since if $\sigma\in {\mathcal F}_{\rm ext}$,  $\int_{\sigma} {\bm u}_K ({\bm x})\cdot {\bm n}_{K,\sigma} {\rm d} s({\bm x}) = 0$, and if ${\mathcal M}_\sigma = \{K,L\}$, then, by definition of $V_h$ from $\widehat{V}_h$, there holds
\[
 \int_{\sigma} {\bm u}_K ({\bm x})\cdot {\bm n}_{K,\sigma} {\rm d} s({\bm x}) + \int_{\sigma} {\bm u}_L ({\bm x})\cdot {\bm n}_{L,\sigma} {\rm d} s({\bm x})=0.
\]
The main advantage of Definition \eqref{eq:defbh} for $b_h({\bm u},{\bm v})$ is the following: for a given ${\bm V}_{\sigma}$, it only involves values ${\bm U}_{\sigma'}$ with $\sigma'\in {\mathcal F}_K$, which means that, using a Newton-Raphson method, the stencil of the Jacobian matrix issued from the trilinear term is block-diagonal, similarly to the diffusion terms (note that it leads to cross dependencies between all the components of the velocities).

\end{enumerate}

\section{Study of the linear systems}\label{sec:matrix}

\subsection{The coupled velocity-pressure linear system}

Let us now detail the construction of the linear system which is directly issued from  \eqref{eq:pweakdis} in the case where $b_h =0$ or issued from the Newton method applied to  \eqref{eq:pweakdis} if $b_h\neq 0$. For any finite set $E$, we denote by $\#E$ its cardinal.

\medskip

This system of linear equations is obtained, first selecting ${\bm v}$ in  the first equation of \eqref{eq:pweakdis} with one component equal to 1 and all the other ones equal to 0, then  selecting  $q$ in  the second equation of \eqref{eq:pweakdis} with  one component equal to 1 and all the other ones equal to 0. 
Letting $U = ((U_{i,\sigma})_{i=1,\ldots,d,\sigma\in {\mathcal F}_{\rm int}}$ and  $P = (P_K)_{K\in {\mathcal M}\setminus\{K_0\}}$, the linear system reads
\begin{equation}\label{eq:syslin}
 \begin{bmatrix} A & D^t \\ D & 0\end{bmatrix} \begin{bmatrix} U\\ P\end{bmatrix}=  \begin{bmatrix} R\\ 0\end{bmatrix}.
\end{equation}
Let us detail the construction of the matrices $A$ and $D$, and of the right-hand side $R$.

\medskip

For any $K\in {\mathcal M}$, we first define the elementary assembly matrix $S_K$, whose side is equal to $s_K :=\#{\mathcal F}_{K,\rm int}$ (recall that this side is equal to $3$ in 2D and $4$ in 3D for any interior element $K$), by
\[
 (S_K)_{\sigma,\sigma'} =  \int_K \Big(\mu  \Pi_h \varphi_\sigma\Pi_h \varphi_{\sigma'} +\nu \nabla\varphi_\sigma\cdot\nabla\varphi_{\sigma'}\Big){\rm d}{\bm x}.
\]
Note that the matrix $S_K$ is symmetric positive definite if $\mu>0$ (transient problems) and only symmetric positive if $\mu=0$.

\medskip

We now define the elementary assembly matrix $A_K$, whose side is equal to $ds_K$, such that, if $b_h = 0$,
\[
 (A_K)_{i,\sigma,j,\sigma'} = \begin{cases}(S_K)_{\sigma,\sigma'} & \hbox{ if }i=j\\ 0 & \hbox{ otherwise }\end{cases}
\]
If $b_h\neq 0$, this matrix is completed with the derivatives of the convection term with respect to the local velocity unknowns $i=1,\ldots,d$ and $\sigma\in {\mathcal F}_{K,\rm int}$.

\medskip

We then define, for any element $K$ of the mesh, the rectangular matrix $H_K$ with $d \#{\mathcal F}_{\rm int}$ lines and $ds_K$ columns, such that, at the column associated to the local velocity unknown $(i,\sigma)\in\{1,\ldots,d\}\times {\mathcal F}_{K,\rm int}$, all the components are null except the one that is at the line associated to the global unknown $U_{i,\sigma}$.

Then the matrix $A$ in \eqref{eq:syslin} is obtained by assembling the elementary matrices, as follows:
\[
 A = \sum_{K\in{\mathcal M}} H_K A_K H_K^t.
\]
For the line of $A$ associated to the global unknown $U_{i,\sigma}$, non-zero terms may occur at the columns associated to the global unknown $U_{j,\sigma'}$ such that there exists $K\in{\mathcal M}$ with $\sigma,\sigma'\in {\mathcal F}_{K,\rm int}$. If $b_h=0$, the matrix $A$ is symmetric positive definite; its inverse is a full matrix, so one cannot solve the linear system by eliminating the velocity unknowns.

\medskip

We define the matrix $D_K$ with $ds_K$ lines and one column (it is then assimilated to a vector), letting for $i\in\{1,\ldots,d\}$ and $\sigma\in {\mathcal F}_{\rm int}$,
\[
 (D_K)_{i,\sigma} = - {\bm a}_{K,\sigma}^{(i)}.
\]
We then define the rectangular matrix $F_K$, with $\#{\mathcal M} - 1$ lines  and 1 column, by 0 everywhere, except 1 at the line corresponding to the global unknown $P_K$, for $K\in {\mathcal M}\setminus\{K_0\}$. Then the matrix $D$ in \eqref{eq:syslin} is defined by
\[
 D = \sum_{K\in{\mathcal M}\setminus \{{K_0}\}} F_K D_K^t H_K^t.
\]

Finally, for any $K\in{\mathcal M}$, let $R_K$ be the elementary right-hand-side issued from \eqref{eq:defrhs}, under the form of a vector with  $ds_K$ components, defined in the case where $b_h=0$,  for all $i\in\{1,\ldots,d\}$ and $\sigma\in {\mathcal F}_{K,\rm int}$ by
\[
 (R_K)_{i,\sigma} =  {\bm a}_{K,\sigma}^{(i)} \Big(\frac 1 {|K|}\int_K\overline{\bm f}{\rm d}{\bm x}\Big)\cdot ({\bm x}_\sigma - {\bm x}_K).
\]
In the case where $b_h\neq 0$, $R_K$ is completed by the non-linear terms issued from the Newton method.
Then the assembled right hand side in \eqref{eq:syslin} is given by
\[
 R = \sum_{K\in{\mathcal M}}H_K R_K.
\]
As recalled in the introduction, the resolution of \eqref{eq:syslin} is then a difficult problem for large meshes. Direct methods can no longer be used, and iterative methods must be based on efficient preconditioners.

\subsection{Hybridisation of the linear system}\label{sec:matstokes}

We construct in this section a linear system, whose the solution directly provides that of \eqref{eq:syslin}, and which can be solved in some cases (see the numerical examples) by cheaper methods. As recalled in the introduction of this paper, the method used for constructing this linear system follows the hybridisation method used in \cite{chen1996equi, duran, vohralik}.

To this purpose, we introduce, for any  $K\in{\mathcal M}$, two diagonal matrices $E_K$ and $C_K$ with the same side  $ds_K$, satisfying the following properties:
\begin{equation}
 (E_K)_{i,\sigma,i,\sigma} + (E_L)_{i,\sigma,i,\sigma}= 0 \hbox{ for all }i=1,\ldots,d,
\label{eq:prope}\end{equation}
and
\begin{equation}
 (C_K)_{i,\sigma,i,\sigma} + (C_L)_{i,\sigma,i,\sigma}= 0 \hbox{ for all }i=1,\ldots,d,
\label{eq:propc}\end{equation}
in the case where ${\mathcal M}_{\sigma} = \{K,L\}$. The matrix $C_k$ is meant to be invertible (in practice, we let the diagonal terms of $C_K$ be equal to $\pm 1$), whereas, if $\mu>0$, the choice $E_K = 0$ can be done.

\medskip

We consider a global vector $\widehat U_{K,i,\sigma}$, associated to the component $i\in\{1,\ldots,d\}$ of the velocity defined at the face $\sigma\in {\mathcal F}_{K,\rm int}$ of $K\in {\mathcal M}$. The number of components of this vector is equal to $\sum_{L\in {\mathcal M}} d \ s_L$; this number is equal to $2 d \#{\mathcal F}_{\rm int}$ since any velocity unknown appears twice at any interior face.

We then define, for any element $K$ of the mesh, in a similar way to the matrix $H_K$, the rectangular matrix $\widehat H_K$ with $\sum_{L\in {\mathcal M}} d \ s_L$ lines and $ds_K$ columns, such that, at the column associated to the local velocity unknown $(i,\sigma)\in\{1,\ldots,d\}\times {\mathcal F}_{K,\rm int}$, all the components are null except the one that is at the line associated to the global unknown $\widehat U_{K,i,\sigma}$.

\medskip

Let us define the following matrices, using the matrices $E_K, C_K, \widehat H_K$ defined in this section and the matrices $A_K, H_K, D_K, F_K$ defined in the previous section:
\[
 \widehat A_K = A_K+E_K\hbox{ and } \widehat A = \sum_{K\in{\mathcal M}} \widehat H_K \widehat A_K   \widehat H_K^t,
\]
\[
 \widehat D = \sum_{K\in{\mathcal M}\setminus \{{K_0}\}} F_K D_K^t  \widehat H_K^t,
\]
\[
\widehat  C = \sum_{K\in{\mathcal M}} H_K   C_K  \widehat H_K^t,
\]
and the following right-hand side, using the right-hand sides $R_K$ defined in the previous section:
\[
\widehat R = \sum_{K\in{\mathcal M}}\widehat H_K R_K.
\]
We consider the following unknown
\begin{itemize}
\item $\widehat{U} = (\widehat U_{K,i,\sigma})$ for $K\in {\mathcal M}$, component $i\in\{1,\ldots,d\}$ and $\sigma\in {\mathcal F}_{K,\rm int}$,
\item $\widehat{P} = (\widehat P_K)$ for $K\in {\mathcal M}\setminus\{K_0\}$,
\item $\widehat{W} = (\widehat W_{i,\sigma})$ for $i\in\{1,\ldots,d\}$ and $\sigma\in {\mathcal F}_{\rm int}$,
\end{itemize}
solution to the following linear system
\begin{equation}\label{eq:syslinhyb}
 \begin{bmatrix} \widehat{A} & \widehat{D}^t & \widehat{C}^t\\ \widehat{D} & 0 & 0\\ \widehat{C} & 0 & 0\end{bmatrix} \begin{bmatrix} \widehat{U}\\ \widehat{P} \\ \widehat{W}\end{bmatrix}=  \begin{bmatrix} \widehat{R}\\ 0\\ 0\end{bmatrix}.
\end{equation}
In the preceding linear system, the equations $\widehat{A} \widehat{U}+ \widehat{D}^t \widehat{P}  + \widehat{C}^t \widehat{W}=\widehat{R}$ can be seen as the splitting of the equations $(AU+ D^t  P)_{\sigma} =R_{\sigma}$ with $P=\widehat{P}$, which hold for all $\sigma\in {\mathcal F}_{\rm int}$, into two equations, one for $K,\sigma$ and the other one for $L,\sigma$ when ${\mathcal M}_{\sigma} = \{K,L\}$, thanks to the introduction of an additional unknown $\widehat{W}_{\sigma}$. 
The velocity unknowns are also doubled, and the equality between the doubled unknowns is ensured by the relation  $\widehat{C}\widehat{U} = 0$.

The next paragraphs are providing details on the following points (among others):
the elimination of $\widehat{W}$ is done by addition of these two equations (owing to \eqref{eq:propc}), and then one recovers  $(AU+ D^t  P)_{\sigma} =R_{\sigma}$ (owing to \eqref{eq:prope}); the system \eqref{eq:syslinhyb} is well-posed, and it is possible, under appropriate choices of the matrices $E_K$, to eliminate $\widehat{U}$ and $\widehat{P}$ in \eqref{eq:syslinhyb}, in order to obtain a linear system only on $\widehat{W}$, with the same stencil as the matrix $A$, and which is symmetric positive definite in some situations.

\medskip

Indeed, the following properties hold.

\begin{enumerate}
\item  {\bf Block diagonal property of $\widehat H_K$ and $\widehat A$.}

 We have the property, for all $K,L\in  {\mathcal M}$,
\begin{equation}
\widehat H_K^t \widehat H_L =\begin{cases}{\rm Id}_K & \hbox{ if }K=L\\ 0 & \hbox{ otherwise. }\end{cases}
\label{eq:proph}\end{equation}

Moreover, the matrix  $\widehat A$ has the blocks $\widehat A_K$ on the diagonal and is null elsewhere. In the case where all the matrices $(\widehat A_K)_{K\in{\mathcal M}}$ are invertible, there holds
\[
 \widehat{A}^{-1} = \sum_{K\in{\mathcal M}} \widehat H_K \widehat A_K^{-1}   \widehat H_K^t.
\]
This leads to a cheap computation of $\widehat{A}^{-1}$ and fully scalable.

\item {\bf Recovery of the solution to  \eqref{eq:syslin}.}

Any solution $(\widehat{U}, \widehat{P} ,\widehat{W})$  of \eqref{eq:syslinhyb} must satisfy
\[
\widehat{C} \widehat{U}=  \sum_{K\in{\mathcal M}} H_K   C_K  \widehat H_K^t\widehat{U} = 0.
\]
For any $i\in\{1,\ldots,d\}$ and $\sigma\in {\mathcal F}_{K,\rm int}$ with ${\mathcal M}_{\sigma} = \{K,L\}$, this means that
\[
 (C_K)_{i,\sigma,i,\sigma} \widehat U_{K,i,\sigma}+ (C_L)_{i,\sigma,i,\sigma}\widehat U_{L,i,\sigma}= 0,
\]
which, together with \eqref{eq:propc} and the invertibility of $C_K$ and $C_L$, provides
\begin{equation}
\widehat U_{K,i,\sigma} =\widehat  U_{L,i,\sigma} :=U_{i,\sigma},
\label{eq:ukl}\end{equation}
denoting by $U_{i,\sigma}$ this common value. Introducing the vector $U = (U_{i,\sigma})_{i=1,\ldots,d,\ \sigma\in {\mathcal F}_{\rm int}}$, we then have
\begin{equation}
\widehat H_K^t \widehat{U} = H_K^t U \hbox{ for all }K\in {\mathcal M}.
\label{eq:propu}\end{equation}
We now multiply by the left the equality
$\widehat{A} \widehat{U}+ \widehat{D}^t  \widehat{P}  + \widehat{C}^t \widehat{W}=\widehat{R}$ by the matrix $J$ which is the matricial translation of the addition of the two equations $K,i,\sigma$ and $L,i,\sigma$ for ${\mathcal M}_{\sigma} = \{K,L\}$. This matrix $J$, which has $d \#{\mathcal F}_{\rm int}$ lines and $2 d \#{\mathcal F}_{\rm int}$ columns, is defined by
\[
 J = \sum_{K\in{\mathcal M}} H_K\widehat H_K^t.
\]
On each line of $J$, all the components are null except two of them, equal to $1$, which enables the addition of pairs of lines.
We then obtain
\[
 J\widehat{A} \widehat{U}+ J\widehat{D}^t  \widehat{P}  + J\widehat{C}^t \widehat{W}=J\widehat{R}.
\]
We then remark that, accounting for \eqref{eq:proph},
\[
 J\widehat{A} \widehat{U} = \sum_{K\in{\mathcal M}} H_K\widehat H_K^t \sum_{L\in{\mathcal M}} \widehat H_L (A_L+E_L)  \widehat H_L^t \widehat{U} = \sum_{K\in{\mathcal M}} H_K (A_K+E_K)  \widehat H_K^t \widehat{U}.
\]
We apply \eqref{eq:propu}, thus obtaining
\[
 J\widehat{A} \widehat{U} = \sum_{K\in{\mathcal M}} H_K (A_K+E_K)  H_K^t U.
\]
Let us now observe that the matrix $\sum_{K\in{\mathcal M}} H_K E_K  H_K^t$ vanishes applying \eqref{eq:prope}. We then get
\[
  J\widehat{A} \widehat{U} = \sum_{K\in{\mathcal M}} H_K A_K  H_K^t U = A U.
\]
We now compute, again accounting for \eqref{eq:proph},
\[
 J\widehat{D}^t \widehat{ P}  = \sum_{K\in{\mathcal M}} H_K\widehat H_K^t \sum_{L\in{\mathcal M}\setminus \{{K_0}\}} \widehat H_K  D_L^t F_L^t \widehat{P}  =  \sum_{K\in{\mathcal M}\setminus \{{K_0}\}} H_K  D_K F_K^t \widehat{P}  = D^t \widehat{ P} ,
\]
\[
 J\widehat{R}=\sum_{K\in{\mathcal M}} H_K\widehat H_K^t\sum_{L\in{\mathcal M}}\widehat H_L R_L =  \sum_{K\in{\mathcal M}}H_K R_K = R.
\]
The matrix $J\widehat{C}^t$ satisfies
\[
 J\widehat{C}^t = \sum_{K\in{\mathcal M}} H_K\widehat H_K^t \sum_{L\in{\mathcal M}}  \widehat  H_L^t   C_L H_L =  \sum_{K\in{\mathcal M}} H_KC_K H_K,
\]
which vanishes owing to  \eqref{eq:propc}.
So we get
\[
 A U +  D^t  \widehat{P}  = R.
\]

Turning to the equation $\widehat{D} \widehat{U} = 0$, we get
\[
\widehat{D} \widehat{U} =  \sum_{K\in{\mathcal M}\setminus \{{K_0}\}} F_K D_K^t  \widehat H_K^t  \widehat{U}=  \sum_{K\in{\mathcal M}\setminus \{{K_0}\}} F_K D_K^t U = D U = 0,
\]
applying \eqref{eq:propu}.
 So we conclude that $(U, \widehat{P} )$ is solution to  \eqref{eq:syslin}. Since this latter system is invertible, we get $\widehat{P} = P$.

 \item{\bf Invertibility of  \eqref{eq:syslinhyb}.}
 
 The invertibility of the linear system is proved, if one assumes that the right-hand side is null, this implies that the solution is null too. This is done by assuming that, in  \eqref{eq:syslinhyb}, we let $\widehat R = 0$ (which is obtained if we let $R_K = 0$ for all $K\in\mathcal{M}$). Since this is a particular case of the linear system under study, the conclusions obtained in the preceding paragraphs, that any solution of this linear system is also a solution to  \eqref{eq:syslin}, are remaining true in this case.  Then, the vectors $U$ issued owing to the preceding computations from $\widehat  U$ and $ \widehat{P} $, are solution to  \eqref{eq:syslin} with $R=0$, since $R$ is computed from null $R_K$. We recall that the linear system  \eqref{eq:syslin} is invertible, which implies that $U = 0$ and $\widehat{ P } = 0$. From $U=0$, we deduce by \eqref{eq:ukl} that $\widehat  U = 0$, which proves from
 $\widehat{A} \widehat{U}+ \widehat{D}^t  \widehat{P}  + \widehat{C}^t \widehat{W}=\widehat{R}$ that 
 \[
 \widehat  C^t \widehat W = 0.
 \]
 The preceding relations are equivalent to $(\widehat C_K)_{i,\sigma,i,\sigma}\widehat  W_{i,\sigma} = 0$ and $(\widehat C_L)_{i,\sigma,i,\sigma}\widehat  W_{i,\sigma} = 0$ , for any $\sigma\in \#{\mathcal F}_{\rm int}$ with ${\mathcal M}_\sigma = \{K,L\}$.
which shows that $\widehat W= 0$ (recall that the matrices $\widehat C$ must have a non-zero diagonal).

The linear system \eqref{eq:syslinhyb} is therefore invertible, and its resolution provides the solution to  \eqref{eq:syslin}.

\item {\bf Elimination of $( \widehat{U},P)$.}

Assuming that, for all $K\in \mathcal{M}$, all the eigenvalues of the symmetric matrix $\widehat{A}_K$ are either strictly positive or strictly negative, let us proceed to the elimination of  $\widehat{U}$ and $P$.
We first have
\[
 \widehat{U} = \widehat{A}^{-1}(- \widehat{D}^t  P  - \widehat{C}^t \widehat{W} +\widehat{R}).
\]
This yields
\[
 \widehat{U} = \sum_{K\in{\mathcal M}} \widehat H_K \widehat{A}^{-1}(- \widehat{D}^t  P  - \widehat{C}^t \widehat{W} +\widehat{R}),
\]

Then we have
\[
 \widehat D\widehat{A}^{-1}(- \widehat{D}^t  P  - \widehat{C}^t \widehat{W} +\widehat{R})= 0.
\]
Let us compute the matrix  $B = \widehat D\widehat{A}^{-1}\widehat{D}^t$.
Using the property
\[
\widehat H_K^t \widehat{A}^{-1}\widehat H_L =\begin{cases}\widehat{A}_K^{-1} & \hbox{ if }K=L\\ 0 & \hbox{ otherwise, }\end{cases}
\]
we get
\[
 B = \sum_{K\in{\mathcal M}\setminus \{{K_0}\}} F_K D_K^t  \widehat{A}_K^{-1}D_K F_K^t.
\]
We then get that $B$ is the diagonal matrix with the values $B_K :=D_K^t  \widehat{A}_K^{-1}D_K$ on the diagonal. Letting ${\underline \lambda}_K$ be the smaller absolute value of the eigenvalues of $\widehat{A}_K^{-1}$, we get that
\[
  |B_K| \ge {\underline \lambda}_K \Vert D_K\Vert_2 >0,
\]
since there exists at least one component of $D_K$ which is different from $0$. So the diagonal matrix $B$ is invertible, and we can write
\[
 B^{-1} = \sum_{K\in{\mathcal M}\setminus \{{K_0}\}} \frac 1 {B_K} F_K F_K^t,
 \]
and
\[
 P = B^{-1}\widehat D\widehat{A}^{-1}( - \widehat{C}^t \widehat{W} +\widehat{R}).
\]
We then obtain
\[
 \widehat C\widehat{A}^{-1}(- \widehat{D}^t  P  - \widehat{C}^t \widehat{W} +\widehat{R}) = 0,
\]
which leads, denoting $G =  \widehat C\Big(\widehat{A}^{-1} - \widehat{A}^{-1}\widehat{D}^t  B^{-1}\widehat D\widehat{A}^{-1}\Big)\widehat{C}^t$ and $S = \widehat C\Big(\widehat{A}^{-1} - \widehat{A}^{-1}\widehat{D}^t  B^{-1}\widehat D\widehat{A}^{-1}\Big)\widehat{R}$, to
\[
G \widehat{W} = S.
\]
The matrix $ G$ is then invertible, since this resolution process is equivalent to the initial linear system (under the above assumption on $\widehat{A}_K$).

\item {\bf Stencil of $G$}

Under the same assumption as previously (for all $K\in \mathcal{M}$, all the eigenvalues of the symmetric matrix $\widehat{A}_K$ are either strictly positive or strictly negative), a simple computation using \eqref{eq:proph} and $F_K^t F_L= 1$ if $K=L$ and $0$ otherwise, leads to
\[
 G  = \sum_{K\in{\mathcal M}} H_K G_K H_K^t,
\]
with, for all $K\in {\mathcal M}\setminus{K_0}$,
\[
 G_K = C_K \Big(\widehat{A}_K^{-1} - \frac 1 {B_K}   \widehat{A}_K^{-1} D_K D_K^t\widehat{A}_K^{-1}\Big)C_K,
\]
and 
\[
 G_{K_0} = C_{K_0} \widehat{A}_{K_0}^{-1} C_{K_0}.
\]
This shows that the assembling of $ G$ leads to the same stencil as that of $A =  \sum_{K\in{\mathcal M}} H_K A_K H_K^t$ (in the case where the matrix $A_K$ is full).

\item  {\bf Case  $\mu>0$ and $b_h=0$.}

In the case  $\mu>0$ and $b_h=0$, all the matrices $A_K$ are symmetric positive definite and we let $E_K  = 0$. Let us show that the resulting matrix $G $ is symmetric positive definite. Indeed, for any vector $\widehat{W}$, let us compute
\[
 a = \widehat{W}^t G \widehat{W}.
\]
Denoting by $Z_K = C_K H_K^t\widehat{W}$, and defining the scalar product $\langle X,Y\rangle_K = X^t\widehat{A}_K^{-1}Y$,  we get that
\[
 a =  \sum_{K\in{\mathcal M}\setminus\{K_0\}} \Big(\langle Z_K,Z_K\rangle_K  - \frac {(\langle Z_K,D_K\rangle_K)^2} {\langle D_K,D_K\rangle_K}\Big)  + \langle Z_{K_0},Z_{K_0}\rangle_{K_0}.
\]
The Cauchy-Schwarz inequality implying
\[
 (\langle Z_K,D_K\rangle_K)^2\le \langle Z_K,Z_K\rangle_K\langle D_K,D_K\rangle_K,
\]
we get that $a\ge 0$.
Since we proved above that, under a weaker hypothesis, the matrix $G$ is invertible, it is then positive symmetric definite.

\item {\bf Computation of $E_K$ in the case $\mu=0$.}

Different strategies can be used. One of them consists in partitioning $\mathcal{M}$ in $\mathcal{M}_1\cup\mathcal{M}_2$, such $\mathcal{M}_2$ is the set of all the neighbours of all $K\in \mathcal{M}_1$. Then for all $K\in \mathcal{M}_1$, we let $E_K = - \lambda{\rm Id}$ with $\lambda$ larger than all the eigenvalues of $A_K$. Then, for all $K\in \mathcal{M}_2$ and $\sigma \in \mathcal{F}_{K,{\rm int}}$ with $\mathcal{M}_\sigma =\{K,L\}$, if  $L\in  \mathcal{M}_1$ (such a $L$ exists by construction), we set $(E_K)_{\sigma,\sigma} = \lambda$. If  $L\notin  \mathcal{M}_1$, we set $(E_K)_{\sigma,\sigma} = 0$.

Then Property \eqref{eq:prope} holds, as well as the fact that all the matrices $(A_K)_{K\in\mathcal{M}_1}$ are symmetric and have all their eigenvalues strictly negative and all  the matrices $(A_K)_{K\in\mathcal{M}_2}$ are symmetric positive definite.

\end{enumerate}

In conclusion of this section, we can use the following method, called the {\bf hybrid method} for solving \eqref{eq:syslin}: 
\begin{enumerate}
 \item One computes the matrix $G$ and the right-hand side $S$ as defined above (this leads to cheap computations).
 \item One then solves the linear system $G \widehat{W} = S$ by a direct method for the small cases or by an iterative method for the larger ones. Note that, in the case where $\mu>0$ and $b_h=0$, a simple preconditioned conjugate gradient solver may be used, and the side of this linear system is smaller than that of  \eqref{eq:syslin} with a stencil similar to that of $A$, which is a part of the matrix of \eqref{eq:syslin}.
 \item One then recovers $P$ and $U$ by the preceding relations which only leads to cheap and fully scalable computations. 
\end{enumerate}

The numerical section provides a few comparisons of this method with the resolution of \eqref{eq:syslin} by a solver with unknowns $(U,P)$.

\section{Numerical results} \label{sec:num}

\subsection{Numerical convergence of the scheme}

Although the Crouzeix-Raviart scheme \eqref{eq:pweakdis} is highly standard in the transient or steady Stokes case, the implementation for the right hand side through the reconstruction $\widehat\Pi_h$ is not completely classical. Note that, if $\overline{f}$ is a constant vector (which means that the velocity is null and that the gradient of the exact pressure is equal to $\overline{f}$), a standard computation of the right-hand side by the integration of $\overline{f}$ against the Crouzeix-Raviart basis functions provides a significant error on the velocity field. On the contrary, owing to the reconstruction $\widehat\Pi_h$, we obtain a null numerical velocity and the exact pressure field (at the machine precision). 

Let us also observe that the non-linear term \eqref{eq:defbh} introduced by \cite{herlat} is not so well-known, and that it is therefore interesting to check, on the analytical Green-Taylor solution, the numerical convergence of this scheme, independently of the algebraic method used for solving the linear systems.

First letting $d=2$, we assume that the analytical solution is given by $ \overline{f} = 0$,
\begin{equation}\label{eq:gretayvel}
 \overline{\bm{u}}({\bm x},t) = {\rm Re}\left(\begin{array}{c}  -\cos(2 \pi (x_1 + \frac 1 4)) \sin(2 \pi (x_2 + \frac 1 2)) \exp(-8 \pi^2 t) 
 \\
 \sin(2 \pi (x_1 + \frac 1 4)) \cos(2 \pi (x_2 + \frac 1 2)) \exp(-8 \pi^2 t) 
 \end{array}\right)
\end{equation}
 and
 \begin{equation}\label{eq:gretaypre}
\overline{p}({\bm x},t)=-\frac {{\rm Re}^2} 4 \Big(\cos(4\pi(x_1 + \frac 1 4))+\cos(4\pi(x_2 + \frac 1 2))\Big)\exp(-16 \pi^2 t).
 \end{equation}
We then implement the values $\overline{\bm{u}}({\bm y},0)$ as initial numerical value at all the nodes of the mesh ${\bm y}$, and the values  $\overline{\bm{u}}({\bm y}_b,t^{(n)})$ at all the boundary nodes of the mesh ${\bm y}_b$ and at the discrete times $t^{(n)} = n\Delta\!t$. The hybrid method and a direct solver are used for these computations which are not dedicated to observe computing performances.
 Letting ${\rm Re} = 100$ and the final time be equal $0.01$, we find the numerical errors given by Table \ref{tab:tab0} with different meshes and time steps.
       
\begin{table}[!h]
  \centering
\begin{tabular}{c|c||c|c|c|c}
  $\Delta\!t$ & $h$ & errl2U & ratio & errl2P & ratio\\
	\hline
  1.25e-04 & 0.2500 &  0.277E+02  &- & 0.487E+03 &-\\
  3.13e-05&0.1250 &  0.854E+01 &  1.70  & 0.262E+03& 0.89 \\
  7.81e-06&0.0625 &  0.315E+01 &  1.44 &     0.112E+03&1.23    \\
  1.95e-06&0.0312 &  0.918E+00 & 1.78&    0.351E+02 &1.67  \\
  4.88e-07&0.0156 &  0.240E+00 & 1.94  &  0.986E+01 &1.83  \\
  1.22e-07&0.0078 &  0.608E-01 & 1.98& 0.299E+01  &1.72
\end{tabular}
\caption{Numerical errors in the case of the 2D Green-Taylor analytical solution of the Navier-Stokes problem.\label{tab:tab0}}
\end{table}
The meshes are  those labelled from 1 to 6 of the triangular family Mesh1 used in the 2D benchmark \cite{2Dbenchmark}. 
The time step $\Delta t$ and the mesh size $h$ are such that $\Delta t/h^2$ is constant.
 The numerical errors are computed at the nodes for the velocities, and at the centre of gravity of the triangles for the pressures. In Table  \ref{tab:tab0}, the ratios are computed by the formula $\log(E_{i-1}/E_{i})/\log(2)$, where $E_{i}$ is a value taken in the column ``errl2U'' or ``errl2P'' and $E_{i-1}$ is the value immediately above in the table.

\medskip

We observe  in Table \ref{tab:tab0} that the numerical order of convergence tends to $2$ for the velocity errors and the finest meshes, and to a value greater than $1$ for the pressure errors, as it is currently observed by numerical schemes in this case.

We now turn out to a 3d case ($d=3$) with ${\rm Re}=100$ and the final time equal to $0.01$. In order to ensure that the 3D meshes present the same regularity factor, the tetrahedral mesh is obtained by splitting in 6 tetrahedra each cube of a uniform cubic mesh of the test domain. The common side of all the cubes of the cubic mesh have all the same side $h$. 

The first 3D  numerical test concerns a Stokes problem case, where the analytical solution is an extension to the 3D case of the preceding Green-Taylor test. The  first two components of the velocity are given by \eqref{eq:gretayvel} extended for all $x_3\in [0,1]$, the third component is equal to $0$ on the whole domain as well as the  pressure (recall that in the Green-Taylor test, the non-linear term is balanced by the pressure gradient). Using the hybrid method, and a conjugate gradient solver with the ``boomer AMG'' preconditioners, we obtain the results provided by Table \ref{tab:3D_test_Stokes_Green_Taylor}.

\begin{table}[!h]
  \centering
\begin{tabular}{c|c||c|c|c|c}
  $\Delta\!t$ &h & errl2U & ratio & errl2P & ratio\\
	\hline
	1.00e-4 & 1.38e-1 & 3.44 & -& 33& -\\
	2.50e-5 & 6.88e-2 &  0.85 & 2.02 & 17 &0.96\\
	6.25e-6& 3.44e-2 &  0.22&  1.95   & 8.1 &1.07  \\
	 1.56e-6&1.72e-2& 5.4e-2 & 2.03  &     4.0&1.02
\end{tabular}
\caption{Numerical errors in the case of the 3D Green-Taylor analytical solution of the Stokes problem.
\label{tab:3D_test_Stokes_Green_Taylor}}
\end{table}

The convergence orders shown in Table \ref{tab:3D_test_Stokes_Green_Taylor} are similar to those observed in Table  \ref{tab:tab0}. Turning to a 3D Navier-Stokes case, we again consider the extension to the 3D case of the 2D Green-Taylor test.  The  first two components of the velocity are again given by \eqref{eq:gretayvel} for any $x_3\in [0,1]$, the third component is again equal to $0$ on the whole domain, and the pressure is given by  \eqref{eq:gretaypre} for any $x_3\in [0,1]$. 
Again, applying the same method for solving the linear systems as in the previous test case, we obtain the results provided by Table \ref{tab:3D_test_Navier-Stokes_Green_Taylor}.
\begin{table}[!h]
  \centering
\begin{tabular}{c|c||c|c|c|c}
  $\Delta\!t$ &h & errl2U & ratio & errl2P & ratio\\
	\hline
	1.56e-4&6.88e-2&10.7&  -   &297  & - \\
	3.91e-5&3.44e-2 &  4.10 &  1.38 &  143   &    1.05\\
	9.77e-6&1.72e-2&1.28 & 1.68& 48.5   & 1.56 \\
\end{tabular}
\caption{Numerical errors in the case of the 3D Green-Taylor analytical solution of the Navier-Stokes problem.
\label{tab:3D_test_Navier-Stokes_Green_Taylor}}
\end{table}

The convergence orders shown in Table \ref{tab:3D_test_Navier-Stokes_Green_Taylor} show a light loss of convergence order in this case, compared to the ones observed in Table  \ref{tab:3D_test_Stokes_Green_Taylor}, although they give a numerical confirmation of the efficiency of the scheme.

\medskip

These tests validate the use of the Crouzeix-Raviart scheme  \eqref{eq:pweakdis} in association with the trilinear term \eqref{eq:defbh}, in 2D and 3D cases. The remaining part of the numerical section is now devoted to 2D and 3D comparisons of the computing performances for solving the linear systems, with or without the use of the hybrid method, in association with a variety of linear solvers.


\subsection{Comparison of algebraic methods and solvers on the transient Stokes problem}

The aim of this section is to assess the interest of the hybrid method in the case of transient Stokes problems (that means that $\mu = 1/\Delta\!t >0$ and $b_h = 0$). In this case, as seen above, the hybrid method leads to positive symmetric definite linear systems, compared to the non-hybrid method, which only provides symmetric linear systems which are not positive and larger.

\medskip

We performed the computation using a direct sequential solver, the only purpose of these tests being to assess the gain of computing time per time step due to smaller linear systems with the hybrid method compared to the linear systems without the hybrid method. We consider the 3D Green-Taylor Stokes problem, with analytical solution given by \eqref{eq:gretayvel} and $p=0$. The linear systems are solved with a simple Gaussian elimination with natural ordering, the time step is equal to $5. 10^{-4}$ and various meshes are used (see Table \ref{tab:3D_test_Stokes_Green_Taylor_direct_solver}). The decrease in the size of the linear systems leads to a clear diminution in the computing time.

\begin{table}[!h]
  \centering
\begin{tabular}{c||c|c}
  Ncv & not hybrid & hybrid\\
	\hline
	46& 2.9e-3 & 3.7e-3 \\
	384& 1.7e-1 & 1.2e-1 \\
	3062& 2.5e+1 & 1.4e+1 \\
	24576& 3.4e+3 & 2.4e+3 \\
\end{tabular}
\caption{Computation time in seconds per linear system solved by Gaussian elimination in the case of the 3D Green-Taylor analytical solution of the Stokes problem. Ncv denotes the number of tetrahedra.
\label{tab:3D_test_Stokes_Green_Taylor_direct_solver}}
\end{table}

\medskip

We now turn to the evaluation of the possibility to use parallel solvers with or without the hybrid method. All the tests are done using the HYPRE/Euclid library for the solvers and the preconditioners, on a computer with 16 cpus.

\medskip

{\bf Conjugate gradient with algebraic multi grid preconditioners in 2D.}

\medskip

We study the possibility of using the BoomerAMG preconditioners, which is known to provide an optimal speed-up in the case of the linear systems issued from diffusion operators. The numerical choices are the following:
\begin{itemize} 
\item The mesh is ``Mesh1-7''  of the triangular family Mesh1 used in the 2D benchmark \cite{2Dbenchmark} (it corresponds to a mesh size equal to $h=3.90625\cdot 10^{-3}$, which leads to $917\,504$ triangles),
\item Smoother algorithm : Hybrid symmetric Gauss-Seidel or SSOR
\item Parallel coarsening algorithm : one-pass Ruge-Stueben coarsening on each processor, no boundary treatment.
\end{itemize}

We observe that, without hybridisation, non-convergence is observed in all tested cases.

On the contrary, using hybridisation, the convergence of the method is obtained. In Table \ref{tab:AMG_CR2}, we provide the computing times needed for the resolution of one linear system (in this transient Stokes problem with constant time step, all the linear systems have the same matrix) for two different values of the time step.

\begin{table}[!htbp]
\centering
\begin{minipage}{\textwidth}
\centering
	\begin{tabular}{|c||c|c||c|c|}
		\hline
	proc & \begin{tabular}{c}$\Delta\!t =0.0001$\\time/iter (s) \end{tabular} & speed-up & \begin{tabular}{c}$\Delta\!t=0.0512$\\time/iter (s) \end{tabular}& speed-up \\
	\hline
1	&163 &-&656&-\\
2	&90&1.81&331&1.98\\
4	&	47&1.91&166&1.99\\
8	&25&1.88&89&1.87\\
16&19&1.32&65&1.37\\
\hline
	\end{tabular}
	\caption{Computation time with hybridisation, using conjugate gradient with boomerAMG}
	\label{tab:AMG_CR2}
\end{minipage}

\end{table}

\medskip

{\bf CG, BCGS et GMRES with ILU in 2D.}

\medskip

We now consider the case where we use different linear solvers (we use ``Mesh1-7'' with $\Delta\!t =0.0001$):

\begin{itemize}
 \item CG : preconditioned conjugate gradient,
 \item BCGS : Bi-conjugate gradient with stabilization,
 \item GMRES,
\end{itemize}

with the Euclid/ ILU preconditioners. Recall that the efficiency of ILU is mainly lost in the case of multi-processor computations, but that it remains in any case much greater that that of boomer AMG. A parameter of ILU is the filling degree (from 1 to 4 in our tests). 

We again observe that no convergence is obtained using conjugate gradient without hybridisation. We show in Table \ref{tab:CG_CR2_new} the results obtained using conjugate gradient with hybridisation. These results show a lower speed-up compared to the use of boomer AMG, but better absolute performances. Let us finally observe that no results were obtained with increasing the filling degree of the ILU method with more than one processor.

\begin{table}[!htbp]
\centering
	\begin{tabular}{|c|c|c|}
		\hline
	proc & ILU & time/iter (s) \\
	\hline
1&	1	&25 \\
1	&2&19 	\\
1	&3&	22 \\
1	&4	&18 \\
\hline
2	&1	&18 \\
\hline
4	&1&	11 \\
\hline
8	&1	&11 \\
\hline
16	&1	&8.4 \\
\hline
	\end{tabular}
	\caption{Computation time with hybridisation, using conjugate gradient with ILU\label{tab:CG_CR2_new}}
\end{table}

We also used the BCGS and GMRES methods without hybridisation. We then get no result with more that 2 processors, the best performance being 23 s per iteration with 4th degree of ILU, BCGS and 1 processor.

\medskip

{\bf Numerical results in 3D}

\medskip

We only obtained numerical results using ILU preconditioners and only one processor.

\medskip

In these conditions, the results without hybridisation with BCGS were better than those with conjugate gradient, whatever be the degree of filling of the ILU method: for example, using degree 2 and BCGS, the time per iteration without hybridisation is 376 s with Ncv = 1\,572\,864 and $\Delta\!t =0.0005$ for the Green-Taylor problem in Stokes conditions, where it is equal to 498 s with conjugate gradient and hybridisation. Additional tests seem to be necessary for improving this comparison.

\subsection{Comparison of linear solvers on the steady lid driven cavity test in 2D}

This test is dedicated to the comparison of the efficiency of the different algebraic solvers in the case of the steady lid driven cavity with ${\rm Re}=1000$, in 2 space dimensions. We again consider the mesh named ``Mesh1-7''  of the triangular family Mesh1 used in the 2D benchmark \cite{2Dbenchmark} (it corresponds to a mesh size equal to $h=3.90625\cdot 10^{-3}$ and Ncv =  $917\,504$ ).

The non-linear system provided by the scheme is approximated by the Newton method. Since the resulting linear systems are no longer symmetric positive, we cannot use the conjugate gradient method; we use the GMRES method  with a convergence threshold equal to $10^{-11}$ in association with an ILU preconditioners with filling degree 2 to 8. This preconditioners has been shown in several tests to provide a sufficient efficiency, letting the filling degree increase \cite{chenier}. Unfortunately, this efficiency falls down on parallel architectures, so this test is only considered with one processor.

\medskip

 In order to assess the additional difficulty issued from the non-linear terms, we first consider the Stokes  problem (in this case, only one Newton iteration is needed, and the linear system is in fact symmetric, but not positive). 

The numerical results presented in Table \ref{tab:tab1} show that the computation time is largely lower with the hybrid method, compared to the results without hybridisation, and that the comparison shows higher contrasts with low filling degree. 

This observation remains true in the Navier-Stokes case. 
To compare the two methods in the Navier-Stokes case, the GMRES threshold has to be reduced to $10^{-8}$ to ensure the convergence of the linear solver when the non-hybrid approach is employed.
The convergence threshold required for the non-linear iterations is equal to $5.10^{-7}$.
In this case and starting from a fluid flow at rest, 10 and 9 Newton iterations are needed respectively  for the scheme with and without hybridisation. 
The results of Table \ref{tab:tab2} show that, despite one additional Newton iteration, the hybrid method converges about twice quicker than the standard approach.
 

\begin{table}[!h]
  \centering
	\begin{minipage}{.495\textwidth}
	  \caption{Stokes lid driven cavity}
		  \centering
    \begin{tabular}{|c|c|c|}
    \hline
    iLU & Hybrid. & time (s)  \\
    \hline
    \multirow{2}{*}{2} &  without      & 2025   \\
          &  with     &  547.9  \\
    \hline
    \multirow{2}{*}{3} &  without     & 747.9 \\
          &  with      & 477.4\\
    \hline
    \multirow{2}{*}{4} &  without     & 391.4\\
          &  with      & 231.6\\
    \hline
		    \multirow{2}{*}{5} &  without     & 211.1\\
          &  with      & 152.9\\
    \hline
    \multirow{2}{*}{6} &  without     & 163.5\\
          &  with      &  89.13 \\
    \hline
    \multirow{2}{*}{7} &  without     & 128.8\\
          &  with      & 65.96\\
    \hline
    \multirow{2}{*}{8} &  without     & 93.04\\
          &  with      & 69.76\\
    \hline
    \end{tabular}%
  \label{tab:tab1}%
	\end{minipage}\hfill	\begin{minipage}{.495\textwidth}
	  \caption{Navier-Stokes lid driven cavity}
      \centering
\begin{tabular}{|c|c|c|}
    \hline
    iLU & Hybrid. & time (s)  \\
    \hline
    \multirow{2}{*}{2} &  without      & 24120	 \\
          &  with      & 9531 
 \\
    \hline
    \multirow{2}{*}{3} &  without     & 9517 
 \\
          &  with     & 4446
 \\
    \hline
    \multirow{2}{*}{4} &  without    & 5420
 \\
          &  with     & 3169 
 \\
    \hline
    \multirow{2}{*}{5} &  without    & 3244 
 \\
          &  with     & 1696 
 \\
		    \hline
    \multirow{2}{*}{6} &  without    & 2332 
 \\
          &  with     & 1028 
 \\
    \hline
    \multirow{2}{*}{7} &  without    & 1703
 \\
          &  with     & 768.0 
 \\
    \hline
    \multirow{2}{*}{8} &  without    & 1398 
 \\
          &  with     & 717.1 
 \\
    \hline
    \end{tabular}%
  \label{tab:tab2}%
	\end{minipage}
	
\end{table}%


\section{Conclusions}

In 2D and on different test cases of the Crouzeix-Raviart scheme, the numerical results show an advantage for using the hybridisation method for solving the coupled linear systems issued from the Newton-Raphson method or from the Stokes problem.

\medskip

In particular, the hybridisation method allow the use of conjugate gradient solvers.

\medskip

Additional tests must be done in 3D in order to assess the influence of hybridisation on solvers performances.
\bibliography{crsolveur}{}
\bibliographystyle{plain}

\end{document}